\documentclass{article}
\usepackage{amssymb,latexsym,amsmath,amsthm,mathrsfs}
\usepackage{graphicx}
\newcommand{\comment}[1]{}

\begin{document}
\title{Speculations on some characteristic properties of numbers\footnote{Presented to the
St. Petersburg Academy on October 9, 1775.
Originally published as
{\em Speculationes circa quasdam insignes proprietates numerorum},
Acta Academiae Scientarum Imperialis Petropolitinae \textbf{4} (1784),
18--30. E564 in the Enestr{\"o}m index.
Translated from the Latin by Jordan Bell, School of Mathematics and
Statistics, Carleton University, Ottawa, Canada.
Email: jbell3@connect.carleton.ca}}
\author{Leonhard Euler}
\date{}
\maketitle

\S 1. There is no doubt that the multitude of all the different fractions,
which can be constituted between the terms $0$ and $1$, is infinite;
whence, because the multitude of all the numbers together is also infinite,
it is apparent for the multitude of all fractions to be still infinitely greater
than this; for, between two numbers, differing by unity, innumerable different fractions may be admitted. 
Here it is taken that the denominators of the fractions can be increased up to
infinity: and if a term is picked which the numerators are not allowed to
exceed, then certainly the number of fractions, which can be constituted
between the terms $0$ and $1$, will be determinate. 
But, as there is some limit, however large this number will be which is taken for the denominators, 
at first sight this question does
not seem that difficult; truly though
if we carefully consider the matter, so many difficulties occur
that a perfect solution of this question seems hardly possible to
hope for here.

\S 2. Now because the fractions, which we are inquiring into here,
must all be different from each other, from any particular denominator
no other fractions may be formed unless not only the numerators
of them are less than the denominator, but also they are prime to it, as otherwise
they could be reduced to a simpler form. Thus as the fraction 
$\frac{15}{24}$ may be reduced to $\frac{5}{8}$,
this fraction cannot be counted with the denominator $=24$,
since it has already been counted with the denominator $8$.
The whole matter is thus reduced to, for any particular denominator,
which may $=D$, assigning the multitude of numbers less than it and which 
have no common divisor with it,
and of course these can be taken as the numerators just for one particular
denominator.
Thus for the denominator $24$ no other numbers are admitted as numerators
besides $1,5,7,11,13,17,19,23$, the multitude of which is $8$, 
and the ratio of this depends on the composition of the number $24$.
For if the denominator $D$ were a prime number, then certainly
all the numbers less than it, the multitude of which is $D-1$, serve
as suitable numerators. Namely the more divisors the denominator
$D$ has, the more greatly is restricted the multitude of numerators.

\S 3. The question arises here: for any given number $D$, to assign
the multitude of numbers which are less than it and also prime to it.
So that this can be presented more easily, let the character $\pi D$ denote
that multitude of numbers which are less than $D$ and which have no common
divisor with it. And indeed it is clear first that if $D$ were a prime
number, then $\pi=D-1$. As we have previously examined the composite
numbers, we shall tabulate the values of this character $\pi D$ for all
numbers not greater than one hundred:

\[
\begin{array}{l|l|l|l|l}
\pi 1=0&\pi 21=12&\pi 41=40&\pi 61=60&\pi 81=54\\
\pi 2=1&\pi 22=10&\pi 42=12&\pi 62=30&\pi 82=40\\
\pi 3=2&\pi 23=22&\pi 43=42&\pi 63=36&\pi 83=82\\
\pi 4=2&\pi 24=8&\pi 44=20&\pi 64=32&\pi 84=24\\
\pi 5=4&\pi 25=20&\pi 45=24&\pi 65=48&\pi 85=64\\
\pi 6=2&\pi 26=12&\pi 46=22&\pi 66=20&\pi 86=42\\
\pi 7=6&\pi 27=18&\pi 47=46&\pi 67=66&\pi 87=56\\
\pi 8=4&\pi 28=12&\pi 48=16&\pi 68=32&\pi 88=40\\
\pi 9=6&\pi 29=28&\pi 49=42&\pi 69=44&\pi 89=88\\
\pi 10=4&\pi 30=8&\pi 50=20&\pi 70=24&\pi 90=24\\
\pi 11=10&\pi 31=30&\pi 51=32&\pi 71=70&\pi 91=72\\
\pi 12=4&\pi 32=16&\pi 52=24&\pi 72=24&\pi 92=44\\
\pi 13=12&\pi 33=20&\pi 53=52&\pi 73=72&\pi 93=60\\
\pi 14=6&\pi 34=16&\pi 54=18&\pi 74=36&\pi 94=46\\
\pi 15=8&\pi 35=24&\pi 55=40&\pi 75=40&\pi 95=72\\
\pi 16=8&\pi 36=12&\pi 56=24&\pi 76=36&\pi 96=32\\
\pi 17=16&\pi 37=36&\pi 57=36&\pi 77=60&\pi 97=96\\
\pi 18=6&\pi 38=18&\pi 58=28&\pi 78=24&\pi 98=42\\
\pi 19=18&\pi 39=24&\pi 59=58&\pi 79=78&\pi 99=60\\
\pi 20=8&\pi 40=16&\pi 60=16&\pi 80=32&\pi 100=40
\end{array}
\] 

\S 4. From this table it is clear that the denominator $2$
provides only one fraction between $0$ and $1$, namely $\frac{1}{2}$;
the denominator $3$ indeed gives $2$; and $4$ gives the two fractions
$\frac{1}{3}$ and $\frac{3}{4}$, and so on. Whence if we want
to continue the denominators not beyond $10$,
the number of all these fractions will be $31$;
if rather we should continue to $20$, the number is $127$;
and proceeding to $30$ the sum of the fractions gives $277$,
as the following table indicates.

\begin{tabular}{c|c}
Max. denom.&Num. fract.\\
10&31\\
20&127\\
30&277\\
40&489\\
50&773\\
60&1101\\
70&1493\\
80&1975\\
90&2489\\
100&3043
\end{tabular}

\S 5. But then clearly if we wanted to admit all fractions
for the denominator $=10$ the maximum number of all the fractions would
be $1+2+3+4+\cdots+9=45$.
Then those which admit reduction should be excluded. Therefore
first the fractions $\frac{2}{4},\frac{3}{6},\frac{4}{8},\frac{5}{10}$,
which are of course $=\frac{1}{2}$, will be excluded; then indeed $\frac{2}{6}$ and
$\frac{3}{9}$,
of course $=\frac{1}{3}$; and likewise
$\frac{4}{6}$ and $\frac{6}{9}$ as they can be made $=\frac{2}{3}$;
also $\frac{2}{8}$ and $\frac{6}{8}$; and finally
$\frac{2}{10}$, $\frac{4}{10}$, $\frac{6}{10}$, $\frac{8}{10}$,
and the number of all these is $14$, and when this is subtracted from
$45$, $31$ remains. Though for most denominators which we will want to admit,
this enumeration would be too extended, nevertheless, let us 
see how it can be carried out.

\S 6. Thus were $D$ the maximum denominator which we admit, the number
of all fractions will plainly be
\[
=1+2+3+4+\cdots+(D-1)=\frac{DD-D}{2}.
\]
Then all the fractions should be excluded whose value is $\frac{1}{2}$,
aside from $\frac{1}{2}$ itself. To this end, $D$ is divided by $2$ and
the quotient, either exactly or approximate and less, shall be $=\alpha$, and it is clear
that the number of fractions which are to be excluded is $=\alpha-1$.
Then for the fractions $\frac{1}{3}$ and $\frac{2}{3}$, let $\frac{D}{3}=\beta$,
with $\beta$ denoting either exactly or approximate and less, and the number
of fractions to be excluded will be $=2(\beta-1)=(\beta-1)\pi:3$. In a similar
way if we put $\frac{D}{4}=\gamma$; then indeed likewise $\frac{D}{5}=\delta$,
$\frac{D}{6}=\epsilon$, etc., namely until the quotients go past unity,
the numbers of fractions to be excluded will then be
\[
(\gamma-1)\pi:4, \quad (\delta-1)\pi:5,\quad (\epsilon-1)\pi:6, \quad \textrm{etc.}
\]
With these removed, the multitude of fractions which are being searched for
which remain will be:
\[
\frac{DD-D}{2}-(\alpha-1)\pi 2-(\beta-1)\pi 3-(\gamma-1)\pi 4-(\delta-1)\pi 5-\textrm{etc.}
\]
So if it were $D=20$, we will have
\[
\begin{split}
&\frac{20}{2}=\alpha=10,\quad \frac{20}{3}=\beta=6,\quad \frac{20}{4}=\gamma=5,\quad \frac{20}{5}=\delta=4,\\
&\frac{20}{6}=\epsilon=3,\quad \frac{20}{7}=\zeta=2,\quad \frac{20}{8}=\eta=2,\\
&\frac{20}{9}=\theta=2,\quad \frac{20}{10}=\iota=2.
\end{split}
\]
Thus here, because $\frac{DD-D}{2}=190$, the multitude of different fractions will be
\[
\begin{split}
&190-9\cdot 1-5\cdot 2-4\cdot 2-3\cdot 4-2\cdot 2-1\cdot 6-1\cdot 4-1\cdot 6-1\cdot 4\\
&=190-63=127,
\end{split}
\]
as we have found above.

\S 7. Therefore all of this investigation rests on this point, 
that, for any given number $D$, 
the value of the character $\pi D$ needs to be found. And indeed we should
first note, as before, that if $D$ is a prime number then it will
be $\pi D=D-1$. Truly if $D$ is a composite number, the determination
of the character $\pi D$ does not turn out to be too arduous; 
namely it will depend on the factors from which the number $D$ is comprised.

\S 8. Thus let $\pi$ denote any prime number, so that it would be
$\pi p=p-1$, and let us search for the value of $\pi p^2$;
it is certainly clear at once that not all the numbers less than it,
the multitude of which is $pp-1$, are prime to it, but 
just those numbers should be excluded which are divisible by $p$, which are:
$p, 2p, 3p, 4p, \textrm{etc.}, (p-1)p$. But the multitude of these
is $p-1$, and when this number is subtracted from $pp-1$, $p(p-1)$ remains, so that
it would be $\pi pp=(p-1)p$. In a similar way, if it were $D=p^3$,
the multitude of numbers less than it is $p^3-1$, whence those
should be excluded which are divisible by $p$, which are
\[
p,2p,3p,4p,\textrm{etc.},p(pp-1),
\]
the multitude of which is $pp-1$, hence it will be
\[
\pi p^3=p^3-1-(pp-1)=p^3-pp=(p-1)pp.
\]
From this it is now easy to see the for any power, it will in general
be $\pi p^n=(p-1)p^{n-1}$.

\S 9. Now let $q$ be another prime number different than $p$, and
let us look for the value of $\pi pq$. First of all therefore, the
multitude of numbers less than $pq$ is $pq-1$, and thus all those
should be excluded from this which are divisible by either $p$ or by $q$.
Indeed the multiples of $p$ will be
\[
p,2p,3p,4p,\ldots,p(q-1),
\]
the multitude of which is $q-1$. In the same way, 
the multitude of the multiples of $q$ will be $p-1$,
and since these would all be different from the first,
the multitude of all numbers to be excluded will be $p+q-2$,
so that it follows here that
\[
\pi pq=pq-1-(p+q-2)=pq-p-q-1=(p-1)(q-1);
\]
from which we obtain this excellent Theorem: If $p$ and $q$ are
different prime numbers, it will always be
\[
\pi:pq=(p-1)(q-1).
\]
This can be further extended in the same way, that if as well $r$ and $s$
were
prime numbers different from the first, it will be
\[
\begin{split}
&\pi pqr=(p-1)(q-1)(r-1) \, \textrm{and}\\
&\pi pqrs=(p-1)(q-1)(r-1)(s-1).
\end{split}
\]

\S 10. Let us now investigate the value of this formula: $\pi ppq$,
where the multitude of all numbers less than $ppq$ is $ppq-1$,
from
which first all the multiples of $p$ should be excluded,
the multitude of which is $pq-1$;
then indeed the multitude of numbers divisible by $q$ is
$pp-1$,
between which however the numbers occur
\[
pq,2pq,3pq,\textrm{etc.},pq(p-1)
\]
which are also divisible by $p$. Because we want to exclude them here,
this should be removed from the final count, so that
this many will remain $pp-1-(p-1)=pp-p$, 
whence we will then obtain
\[
\pi ppq=ppq-1-(pq-1)-(pp-p)=(p-1)(q-1)p.
\]
Like how it is
\[
p(p-1)=\pi pp \quad \textrm{and} \quad q-1=\pi:q,
\]
this theorem can here be obtained: If $p$ and $q$ are different
prime numbers,
then it will be
\[
pi ppq=\pi pp \cdot \pi q=p(p-1)(q-1).
\]

\S 11. In a similar way it is hardly difficult to see that
\[
\pi p^n q=\pi p^n \pi q=p^{n-1}(p-1)(q-1).
\]
For, because the multitude of numbers less than it is $p^n q-1$,
first here all the multiples of $q$ should be excluded, the number of which is
$p^n-1$, and the multitude that will remain is $p^n q-p^n$.
Besides indeed
we should also exclude all the numbers divisible by $p$, the multitude of which is $p^{n-1}q-1$, and $p^nq-p^n-p^{n-1}q+1$ would remain. To this, however,
all the terms divisible by $pq$ should be added, the multitude of which is
$p^{n-1}-1$, from which one gathers
\begin{eqnarray*}
\pi p^nq&=&p^nq-p^n-p^{n-1}q+p^{n-1}\\
&=&p^{n-1}(pq-p-q+1)=p^{n-1}(p-1)(q-1).
\end{eqnarray*}

\S 12. In a not at all dissimilar way, if a number $D$ were a product from
two powers of any different prime numbers $p$ and $q$,
so that it would thus be $D=p^\alpha q^\beta$,
then it will be
\[
\pi p^\alpha q^\beta=p^{\alpha-1}q^\beta (p-1)(q-1);
\]
and then generally, if the letters $p,q,r,s$ denote prime numbers different
from each other, it will be
\[
\pi p^\alpha q^\beta r^\gamma s^\delta=p^{\alpha-1}q^{\beta-1}r^{\gamma-1}s^{\delta-1}(p-1)(q-1)(r-1)(s-1)
\]
from which one realizes that it will also be
\[
\pi p^\alpha q^\beta r^\gamma s^\delta=\pi p^\alpha \cdot \pi q^\beta \cdot \pi r^\gamma \cdot \pi s^\delta.
\]
Because of this, if only the values of the character $\pi D$ were found for
all powers of prime numbers, then it is perfectly clear that from these,
the values of the character $\pi$ of all numbers could be readily
assigned.

\S 13. If, by means of these Theorems, one wants to investigate the values for
arbitrarily large numbers, the goal will be obtained most quickly
if one resolves the given number $D$ into factors which are prime to each other,
either prime numbers or not. For in fact if it were $D=PQRS$ etc. and these
factors $P,Q,R,S$ have no common divisors, then it will always be
\[
\pi PQRS=\pi P\cdot \pi Q\cdot \pi R\cdot \pi S.
\]
Namely if it were $D=PQRS$ etc., and these factors $P,Q,R,S$ have no
common divisors, then it will always be
\[
\pi PQRS=\pi Q\cdot \pi Q\cdot \pi R\cdot \pi S.
\]
Like if this number were proposed: $D=360$, because $360=9\cdot 40$, it will
be $\pi 360=\pi 9\cdot \pi 40=6\cdot 16=96$.

\S 14. But if indeed the progression of these numbers, which were exhibited in the table given above, are considered, which is $0,1,2,2,4,2,6,4,6,4,10,4$, etc.,
one can find no clear order in the terms of it;
yet in the progression of numbers each term of which exhibits the sums
of the divisors of the natural numbers, I did succeed in detecting a 
 characteristic order. 
Thus at least, if from these numbers such a series were formed:
\[
1x^2+2x^3+2x^4+4x^5+2x^{10}+\textrm{etc.}
\]
the general term of which is signified by our method as $x^n \pi n$, one
sees that
the character of it, or even the sum, might be expressed in 
in some way by known quantities, either algebraic or transcendental.
Therefore it is worth the greatest effort to inquire 
into the nature of this progression, since here the science of numbers
can be enriched with a not negligible increase.

\S 15. However from the general form given in \S 12, a much easier rule
can be deduced, by means of which for any given number $N$ the value
of the character $\pi N$ can be assigned, which we shall explain in the following Problem.

\begin{center}
{\Large Problem.}
\end{center}

{\em Given any number $N$ to find the multitude of all numbers less than it
and prime to it.}

\begin{center}
{\Large Solution.}
\end{center}

\S 16. For any number $N$, it can always be represented in such a form as
$N=p^\alpha q^\beta r^\gamma s^\delta$ etc., with $p,q,r,s$ being prime numbers.
We have also found for it then to be
\[
\pi N=p^{\alpha-1}q^{\beta-1}r^{\gamma-1}s^{\delta-1}(p-1)(q-1)(r-1)(s-1).
\]
Then it will therefore be
\[
\frac{\pi N}{N}=\frac{(p-1)(q-1)(r-1)(s-1)}{pqrs},
\]
from which it follows that
\[
\pi N=\frac{N(p-1)(q-1)(r-1)(s-1)}{pqrs};
\]
so that there does not have to be any more work to know the exponents
$\alpha,\beta,\gamma$, but rather it suffices to just investigate all the different
prime numbers $p,q,r,s$ by which the given number $N$ is divisible;
with these known, the multitude of numbers which are less than $N$ and 
also prime to it will be
\[
\pi N=\frac{N(p-1)(q-1)(r-1)(s-1)}{pqrs}.
\]

\S 17. So if, e.g., this number were proposed: $N=9450$, 
the prime numbers which divide this number are $2,3,5,7$;
since it does not admit division by any other, it will therefore be
\[
\pi 9450=\frac{9450\cdot 1\cdot 2\cdot 4\cdot 6}{2\cdot 3\cdot 5\cdot 7}=2160.
\]

\S 18. Thus if ever $N$ has just a single prime divisor $p$, which happens
when either when $N$ is equal to $p$ itself, or some power of it;
then it will therefore always be $\pi N=\frac{N(p-1)}{p}$. Namely
if it were $N=p$, it will be $\pi N=p-1$; and if it were $N=p^n$, then it will
be $\pi N=p^{n-1}(p-1)$, as we have found above. But if however
$N$ admits two prime divisors $p$ and $q$, then it will be 
$\pi N=\frac{N(p-1)(q-1)}{pq}$. Thus if $N$ has no other divisors besides $2$
and $3$, it will  be $\pi N=\frac{N}{3}$. Such numbers up to one hundred are
\[
6,12,18,24,36,48,54,72,96.
\]

\S 19. For let us take the number $N$ to have the prime divisors $p,q$ and $r$,
different from each other, and besides these no others; and because
the multitude of all numbers not greater than it is $=N$, 
and therefore some number will be divisible by $p,q$ and $r$,
where first all shall be excluded which are divisible by $p$,
the multitude of which would be $\frac{N}{p}$,
and with these deleted the multitude of the remaining will be $N-\frac{N}{p}=\frac{N(p-1)}{p}$;
from this now we should exclude all which are divisible by $q$,
the multitude of which would be $\frac{N(p-1)}{pq}$,
and now there will remain $\frac{N(p-1)(q-1)}{pq}$.
Finally now those which are divisible by $r$ should be excluded,
the multitude of which would be $\frac{1}{r}$ part of this number.
With these
deleted, the number of the remaining will be $\frac{N(p-1)(q-1)(r-1)}{pqr}$;
and in this way our rule has been firmly demonstrated.

\S 20. But nevertheless, this rule provides no help to us on
the nature of the progression which the numbers $\pi N$ constitute,
and which is:
\[
\begin{array}{llllllllllll}
1&2&3&4&5&6&7&8&9&10&11&12\\
0&1&2&2&4&2&6&4&6&4&10&4
\end{array}
\]
which is be explored.
Certainly, if we adjoin powers of the indefinite quantity $x$, and we set
\[
s=1x^2+2x^3+2x^4+4x^5+2x^6+6x^7+4x^8+6x^9+4x^{10}+\textrm{etc.}
\]
from it we can form the following series:
\[
\frac{\int sds}{x}=\frac{xx}{2}+\frac{2x^3}{3}+\frac{x^4}{2}+\frac{4x^5}{5}+\frac{x^6}{3}+\frac{6x^7}{7}+\frac{x^8}{2}+\frac{2x^9}{3}+\frac{2x^{10}}{5}+\textrm{etc.}
\]
where all the coefficients are contained in the formula $\frac{(p-1)(q-1)(r-1)}{pqr}$.

\S 21. Now, all those powers of $x$ will have the same coefficient $\frac{1}{2}$
whose exponents admit just one prime divisor, and thus are powers of
two, namely
\[
x^2,x^4,x^8,x^{16},x^{32},x^{64},\textrm{etc.}
\]
Then all the powers whose exponents are ranks of three, which are
$x^3,x^9,x^{27}$, would all have the same coefficient $\frac{2}{3}$. In a similar way $\frac{4}{3}$ will be the common coefficient of the powers
$x^5,x^{25},x^{125}$, etc. And truly $\frac{1}{3}=\frac{1}{2}\cdot \frac{2}{3}$ will be the common coefficient of all the powers whose exponents involve
exactly the two prime numbers $2$ and $3$, which are
$x^6.x^{12},x^{18},x^{24},x^{36}$, etc. And the same kind of thing happens
for the other exponents, which involve either pairs, or triples, or
quadruples of prime numbers.
Moreover when more prime numbers occur in the exponents, 
the series of powers, which enjoy common coefficients, will be more
plentiful.

\S 22. Thus in this order the simplest series are those
whose constitute a geometric progression,
of which type is $x+x^2+x^4+x^8+x^{16}+\textrm{etc.}$, 
but even the sum of this series has still not been able to be found in any way,
or even to be reduced to some integral formula, and at the very least it is
hoped that some certain order can be found in this series in general,
from which at least the following terms could be determined from the preceding;
this rightly would be seen as all the more remarkable,
since
the coefficient of any power $x^n$ can nevertheless be assigned easily.

\end{document}